\documentclass[onecolumn, leqno, 12pt]{article}
\newtheorem{definition}{Definition}
\newtheorem{remark}[definition]{Remark}

\newtheorem{theorem}[definition]{Theorem}
\newtheorem{proposition}[definition]{Proposition}

\newtheorem{example}[definition]{Example}
\usepackage[pdftex]{hyperref}
\hypersetup{pdfnewwindow=true}
\usepackage{amssymb,amsmath}

\begin{document}
\global\def\refname{{\large \bfseries References:}}
\baselineskip 12.5pt
%
%
%
\title{\bf An Elementary and Real Approach to Values of the  
Riemann Zeta Function\footnote{extended version of the contribution to the proceedings of the XXVII International Colloquium on Group Theoretical Methods in Physics, Yerevan, Armenia, August 13-19, 2008}}

\date{}

\author{\hspace*{-10pt}
\begin{minipage}[t]{2.7in} \normalsize 
\centerline{\sc Armen Bagdasaryan \footnote{\textit{Institution of the Russian Academy of Sciences, 
V.A. Trapeznikov Institute for Control Sciences of RAS, Moscow, Russia}; E-mail: \tt \href{mailto: abagdasari@hotmail.com}{abagdasari@hotmail.com}}}
\end{minipage} \kern 0in
%
%
\\ \\ \hspace*{-10pt}
\begin{minipage}[b]{4.0in} \normalsize
\small An elementary approach for computing the values at negative integers of the Riemann zeta function is presented. The approach is based on a new method for ordering the integers and a new method for summation of divergent series. We show that the values of the Riemann zeta function can be computed, without using the theory of analytic continuation and functions of complex variable.
\\ [4mm] {\footnotesize \textit{MSC numbers:}}
\footnotesize 11M06,  
11B68, 
40C15 
\\ [4mm] {\footnotesize \textit{PACS numbers:}}
02.10.De, 
02.30.Lt 
\end{minipage}
\vspace{-10pt}}

\maketitle
\baselineskip 17,5pt
%
%

\section{Introduction}
\label{S1}\vspace{-4pt}

The Riemann zeta function $\zeta(s)$ is one of the most important objects in the study of number theory (if \textit{not} in the whole modern mathematics). It is a classical and well-known result that $\zeta(s)$, originally defined on the half plane $\Re(s)>1$, can be analytically continued to a meromorphic function on the entire complex plane with the only pole at $s=1$, which is a simple pole with residue 1. 

One of the main reasons of interest to $\zeta(s)$ is that the special values of $\zeta(s)$ at integers have been proved or conjectured to have significant arithmetic meanings. For instance, Zagier's \cite{Zag, Zagg} conjecture concerning the relation between $\zeta(n)$ and the $n$-logarithms for $n\geq 2$  and Lichtenbaum's \cite{Li} conjecture connecting $\zeta(n)$ with motivic cohomology.

Although the Riemann zeta function is considered as being primarily relevant to the "purest" of mathematical disciplines, number theory, it also occurs, for instance, in applied statistics (Zipf's law, Zipf-Mandelbrot law), physics, cosmology \cite{cosmo, cosmo1, cosmo2}. For example, its special values $\zeta(3/2)$ is employed in calculating the critical temperature for a Bose-Einstein condensate, $\zeta(4)$ is used in Stefan-Boltzmann law and Wien approximation. The Riemann zeta function appears in models of 
quantum chaos \cite{berry2,berry1,keat,S} and shows up explicitly in the calculation of the Casimir effect \cite{casimir1, casimir2}. 

The method of zeta function regularization \cite{zetareg}, which is based on the analytic continuation of the zeta function in the complex plane, is used as one possible means of regularization of divergent series in quantum field theory \cite{zetaquantfield4,zetaquantfield3,zetaquantfield1,zetaquantfield2}. 

So, the Riemann zeta function is an interesting object for study not only for mathematicians but for physicists as well. 

\section{Zeta Function and its Values at Negative Integers}
\label{S2}\vspace{-4pt}

The zeta function was first introduced by Euler and is defined by
$$
\zeta(s)=\sum_{u=1}^{\infty}\frac{1}{u^s}.
$$
The series is convergent when $s$ is a complex number with $\Re(s)>1$.

In 1859 Riemann defined $\zeta(s)$ for all complex numbers $s$ by analytic continuation. Several techniques permit to extend the domain of definition of the zeta function (the continuation is independent of the technique used because of uniqueness of analytic continuation). One can, for example, consider the zeta alternating series, so-called the Dirichlet eta function
$$
\eta(s)=\sum_{u=1}^{\infty}\frac{(-1)^{u-1}}{u^s},
$$
defining an analytic function for $\Re(s)>0$. When the complex number $s$ satisfies $\Re(s)>1$, we have
$$
\eta(s)=\sum_{u=1}^{\infty}\frac{1}{u^s}-\sum_{u=1}^{\infty}\frac{2}{(2u)^s}=\zeta(s)-\frac{2}{2^s}\zeta(s)
$$
or
$$
\zeta(s)=\frac{\eta(s)}{1-2^{1-s}}, \;\;\;\;\;\;\;\;\; \Re(s)>1.
$$

Since $\eta(s)$ is defined for $\Re(s)>0$, this identity permits to define the zeta function for all complex numbers $s$ with positive real part, 
except for $s=1+2n\pi i/\log 2$ where $n$ is an integer. For $n=0$, we have a pole at $s=1$, and for $n\neq 0$ one can use the derivative of the $\eta$-function since it is known that $\eta$ has zeros at these isolated points on the line $\Re(s)=1$ \cite{landau,noah,son2,widder}.

So, the zeta alternating series is linked with the original series (the zeta function) by the simple relation
\begin{equation}
\widetilde{\zeta}(s)=(1-2^{1-s})\zeta(s)  \label{eq:zetarelation}
\end{equation}
where we have put $\widetilde{\zeta}(s)=\eta(s)$.

Around 1740, Euler \cite{E, E1, E2} discovered a method of calculating the values of the divergent series
$$
1+1+1+1+...=-\frac{1}{2}
$$
$$
1+2+3+4+...=-\frac{1}{12}
$$
$$
1+4+9+16+...=0
$$
$$
1+8+27+64+...=\frac{1}{120}
$$
etc.

In modern terms, these are the values at non-positive integer arguments of the Riemann zeta function $\zeta(s)$, which, as it was said above, defined by the series, absolutely convergent in $\Re(s)>1$
$$
\zeta(s)=1+\frac{1}{2^s}+\frac{1}{3^s}+\frac{1}{4^s}+...
$$
To find these values, Euler considered the alternating series
\begin{equation}
1^m-2^m+3^m-4^m+5^m-6^m+...=\frac{1}{\widetilde{\zeta}(s)}. \label{eq:euleralternate}
\end{equation}

He observed that the value (\ref{eq:euleralternate}) is obtained as a limit of the power series 
$$
1^m-2^mx+3^mx^2-4^mx^3+...
$$
as $x\rightarrow 1$, since, although the series itself converges only for $|x|<1$, it has an expression as a rational function (analytic continuation, as we now put it), finite at $x=1$, which is obtained by a successive application of multiplication by $x$ and differentiation (equivalently, applying the Euler operator $x\cdot d/dx$ successively after once multiplied by $x$) to the geometric series expansion
\begin{equation}
\frac{1}{1+x}=1-x+x^2-x^3+x^4-x^5+... \;\;\;\;\;\;\; (|x|<1) \label{eq:eulerpower}
\end{equation}

For instance, if we substitute $x=1$ in (\ref{eq:eulerpower}), we find formally
$$
\frac{1}{2}=1-1+1-1+1-1+...=\widetilde{\zeta}(0)
$$
and hence, in view of (\ref{eq:zetarelation}), we have $\zeta(0)=-\frac{1}{2}$. More examples are
$$
\frac{1}{(1+x)^2}=1-2x+3x^2-4x^3+5x^4-...
$$
$$
\frac{1-x}{(1+x)^3}=1-2^2x+3^2x^2-4^2x^3+5^2x^4-...
$$
$$
\frac{1-4x+x^2}{(1+x)^4}=1-2^3x+3^3x^2-4^3x^3+5^3x^4-...
$$
which give us
$$
\widetilde{\zeta}(-1)=\frac{1}{4}, \;\;\;\;\; \widetilde{\zeta}(-2)=0, \;\;\;\;\; 
\widetilde{\zeta}(-3)=-\frac{1}{8}
$$
and in turn
$$
\zeta(-1)=-\frac{1}{12}, \;\;\;\;\; \zeta(-2)=0, \;\;\;\;\; 
\zeta(-3)=\frac{1}{120}.
$$

Euler found these values of divergent series, having no notion of the analytic continuation and functions of complex variable.

In modern concepts this method provides no rigorous way for obtaining the values of $\zeta(s)$ at non-positive integers because it is commonly considered that the values of $\zeta(-m)$ should be established as values of the analytically continued function $\zeta(s)$ at $s=-m$.

There are different approaches to evaluation of the values of $\zeta(-m)$. These values can be computed by using the functional equation satisfied by this
function \cite{Ti}, by using the Euler-Maclaurin formula \cite{G,son}, in \cite{rama} the values of $\zeta(-m)$ are derived from a particular series 
involving $s\zeta(s+1)$, another method is presented in \cite{murty}. Other approaches are based on the so-called $q$-extensions ($q$-analogs) of Riemann zeta function \cite{kaneko,kim1, kim2}.

However, in this paper we find in a mathematically rigorous way a new approach to the values of $\zeta(s)$ at non-positive integer points, without addressing the theory of analytic continuation and functions of complex variable.

\section{Definitions and Preliminaries}
\label{S3}\vspace{-4pt}

We shall deal with two series
$$
\widetilde{\zeta}(s)=\sum_{u=1}^{\infty}\frac{(-1)^{u-1}}{u^s} \;\;\;\;\;\;\; 
\mathrm{(Dirichlet\,\,eta\,\,function)}
$$
and
$$
\zeta(s)=\sum_{u=1}^{\infty}\frac{1}{u^s} \;\;\;\;\;\;\; 
\mathrm{(Riemann\,\,zeta\,\,function)}
$$
and aim to evaluate $\zeta(s)$ and $\widetilde{\zeta}(s)$ at non-positive integers ($\widetilde{\zeta}(-m), \zeta(-m)$) in elementary and real way, without using complex-analytical techniques.

As a starting point we take a new method for ordering the integers \cite{Va, VBT}, which provides very well not only a real and elementary approach to computing the values of $\zeta(s)$ at negative integers but also a potentially new field of research \cite{Ba,VB,VBT}.

In this paper we restrict ourselves to considering the values of $\zeta(s)$ and $\widetilde{\zeta}(s)$ at non-positive integers and leave the issue of computing the values of $\zeta(s)$ and $\widetilde{\zeta}(s)$ at positive integers for our next work \cite{journee}.

To make paper self-contained, we introduce some basic definitions and propositions \cite{VBT} (proofs can be found therein), necessary for the aim of this paper.

To be more specific, let us consider the set of all integer numbers $Z$.\footnote{we denote by $a,b,c$ integer numbers and by $n,m,k$ natural numbers}
\begin{definition}
We say that $a$ precedes $b$, $a, b \in Z$, and write $a\prec b$, if the inequality $\frac{-1}{a}<\frac{-1}{b}$ holds; $a\prec b \Leftrightarrow \frac{-1}{a}<\frac{-1}{b}$ \footnote{assuming by convention $0^{-1}=\infty$}.
\end{definition}

This method of ordering, obviously, gives that any positive integer number (including zero) precedes any negative integer number, and the set $Z$ has zero as the first element and $-1$ as the last element, i. e. we have $Z=[0, 1, 2,...-2, -1]$ \footnote{geometrically, the set $Z$ can be represented as cyclically closed (closed number line)}. In addition, the following two necessary conditions of axioms of order hold:
\begin{enumerate}
\item either $a\prec b$ or $b\prec a$ 
\item if $a\prec b$ and $b\prec c$ then $a\prec c$ 
\end{enumerate}
\begin{definition}
A function $f(x)$, $x\in Z$ , is called \textit{regular} if there exists an elementary function $F(x)$ such that $F(z+1)-F(z)=f(z), \; \forall z\in Z$. The function $F(x)$ is said to be a \textit{generating function} for $f(x)$.
\end{definition}
\begin{remark}
If $F(x)$ is a generating function for $f(x)$, then the function $F(x)+C$, where $C$ is a constant, is also a generating function for $f(x)$. So, any function $F(x)$ which is generating for $f(x)$ can be represented in the form $F(x)+C(x)$, where $C(x)$ is a periodic function with the period $1$. 
\end{remark}

Suppose $f(x)$ is a function of real variable defined on $Z$ and $Z_{a, b}$ is a part of $Z$ such that $Z_{a, b}=[a, b]$ if $a\preceq b$ and $Z_{a, b}=Z\setminus (b, a)$ if $a\succ b$, where $Z\setminus (b, a)=[a, -1]\cup[0, b]$. 
\begin{definition} \label{def:sum}
For any $a, b\in Z$
\begin{equation}
\sum_{u=a}^b{f(u)}=\sum_{u\in Z_{a, b}}{f(u)}\footnote{observing the established order of elements $Z_{a, b}$}.
\end{equation}
\end{definition}

This definition satisfies the condition of generality and has a real sense for any integer values of $a$ and $b$ ($a ^{>}_{<} b$).
The definition \ref{def:sum} extends the classical definition of sum $\sum_{a}^{b}f(n)$ to the case $b<a$.

We introduce for regular functions the following quite natural conditions:
\begin{enumerate}
\item If $S_n=\sum\limits_{u=a}^n{f(u)} \;\; \forall n$, then $\lim\limits_{n\rightarrow \infty}S_n=\sum\limits_{u=a}^\infty {f(u)}$\footnote{$n\rightarrow \infty$ means that $n$ unboundedly increases, without changing the sign}. \label{cond:1}  \\
\item If $S_n=\sum\limits_{u=1}^{n/2}{f(u)} \;\; \forall n$, then $\lim\limits_{n\rightarrow \infty}S_n=\sum\limits_{u=1}^\infty {f(u)}$.\\
\item If $\sum\limits_{u=a}^\infty {f(u)}=S$, then $\sum\limits_{u=a}^\infty {af(u)}=aS, \; a\in R$.\\
\item  If $\sum\limits_{u=a}^\infty {f_1(u)}=S_1$ and $\sum\limits_{u=a}^\infty {f_2(u)}=S_2$, then  $\sum\limits_{u=a}^\infty {(f_1(u)+f_2(u))}= S_1+S_2$.\\
\item For any $a$ and $b$, $a\leq b$: $F(b+1)-F(a)=\sum\limits_{u=a}^b{f(u)}$.\\
\item If $G=[a_1, b_1]\cup[a_2, b_2]$, $[a_1, b_1]\cap[a_2, b_2]=\emptyset$, then \\ $\sum\limits_{u\in G}{f(u)}=\sum\limits_{u=a_1}^{b_1}{f(u)}+\sum\limits_{u=a_2}^{b_2}{f(u)}$. \label{cond:6}  \\
\end{enumerate}

The conditions (\ref{cond:1})-(\ref{cond:6}) define a method of summation of infinite series, which is regular due to (\ref{cond:1}). 

\begin{proposition}
If $f(x)$ is a regular function and $a\in Z$ is a fixed number, then 
\begin{equation}
\sum_{u=a}^{a-1}f(u)=\sum_{u\in Z}f(u)  \label{prop:5}
\end{equation}
\end{proposition}

\begin{proposition}
For any numbers $m$ and $n$ such that $m\prec n$
\begin{equation}
\sum_{u=m}^{n}f(u)=\sum_{u=-n}^{-m}f(-u)  \label{prop:6}
\end{equation}
\end{proposition}

\begin{proposition}
Let $f(x)$ be a regular function and let $a$, $b$, $c$ be any integer numbers such that $b\in Z_{a, c}$. Then 
\begin{equation}
\sum_{u=a}^{c}f(u)=\sum_{u=a}^{b}f(u)+\mathop{{\sum}'}_{u=b+1}^{c}f(u)\;\footnote{the prime on the summation sign means that 
$\mathop{{\sum}'}\limits_{u=b+1}^{c}f(u)=0$ for $b=c$}
\end{equation}
\end{proposition}

\begin{proposition}
Suppose $f(x)$ is a regular function. Then
\begin{equation}
\sum_{u=a}^{a-1}f(u)=0 \;\;\;\;\; \forall{a\in Z}
\end{equation}
or, which is the same in view of (\ref{prop:5})
\begin{equation}
\sum_{u\in Z}f(u)=0 
\end{equation}
\end{proposition}

\begin{proposition}
For any regular function $f(x)$
\begin{equation}
\sum_{u=a}^{b}f(u)=-\sum_{u=b+1}^{a-1}f(u) \;\;\; \forall a, b\in Z  \label{prop:10}
\end{equation}
\end{proposition}

From (\ref{prop:10}), letting $a=0$ and $b=-n$, we have
$$
\sum\limits_{u=0}^{-n}f(u)=-\sum\limits_{u=-n+1}^{-1}f(u)
$$
and using (\ref{prop:6}), we get
$$
\sum\limits_{u=0}^{-n}f(u)=-\sum\limits_{u=1}^{n-1}f(-u)
$$
and
\begin{equation}
\sum\limits_{u=1}^{-n}f(u)=-\sum\limits_{u=0}^{n-1}f(-u) \label{eq:11}
\end{equation}

Using (\ref{eq:11}), we obtain
\begin{theorem}   \label{theor:infsum}
For any even regular function $f(x)$
\begin{equation}
\sum_{u=1}^{\infty}f(u)=-\frac{f(0)}{2}  
\end{equation}
independently on whether the series is convergent or not in a usual sense.
\end{theorem}

\begin{example}
Consider some examples of both convergent and divergent series.
\begin{enumerate}
\item Convergent series 
$$
\sum_{u=1}^{\infty}\frac{1}{4u^2-1}=\frac{1}{2} \;\;\;\;\;\;\;\;\;\;\;\;\;\;\;\;\;\;\;\;\;\;\;\;\;\;\;\;\;\;\;\;\;  \textup{\cite{GR}}
$$

$$
\sum_{u=1}^{\infty}(-1)^{u}\frac{2u^2+1/2}{(2u^2-1/2)^2}=-1 \;\;\;\;\;\;\;\;\;\;\;\;\;\;\;\;\;\;\;\; \textup{}
$$

$$
\sum_{u=1}^{\infty}\frac{(4^u-1)(u-1/2)-1}{2^{u^2+u+1}}=\frac{1}{4} \;\;\;\;\;\;\;\;\;\;\;\;\;\;\;\;\;\;\;\; \textup{}
$$

$$
\sum_{u=1}^{\infty}\frac{(u^2+1/4)\tan(1/2)\cos u-u\sin u}{(4u^2-1)^2}=-\frac{\tan(1/2)}{8} \;\;\;\;\;\;\;\;\;\;\;\;\;\;\; \textup{}
$$
with the generating functions, respectively
$$
F(n)=\frac{-1}{2(2n-1)}, \;\;\; F(n)=\frac{(-1)^{n-1}}{(2n-1)^2}, 
$$
$$
F(n)=\frac{-(n-1/2)}{2^{n^2-n+1}}, \;\;\; F(n)=\frac{\sin(n-1/2)}{8(2n-1)^2\cos(1/2)}
$$

\item Divergent series
$$
1-1+1-1+...=\frac{1}{2} \;\;\;\;\;\;\;\;\;\;\;\;\;\;\;\;\;\;\;\;\;\;\;\;\;\;\;\;\;\;\;\; \textup{\cite{Ha}}
$$

$$
1+1+1+1+...=-\frac{1}{2} \;\;\;\;\;\;\;\;\;\;\;\;\;\;\;\;\;\;\;\;\;\;\;\;\;\;\;\;\; \textup{\cite{Ti}}
$$

$$
1^{2k}+2^{2k}+3^{2k}...=0 \;\;\;\;\; \forall k \;\;\;\;\;\;\;\;\;\;\;\;\;\;\;\;\;\;\;\;\;\;\;\;\;\; \textup{\cite{Ti}}
$$

$$
1^{2k}-2^{2k}+3^{2k}-...=0 \;\;\;\;\; \forall k \;\;\;\;\;\;\;\;\;\;\;\;\;\;\;\;\;\;\;\;\;\; \textup{\cite{Ha}}
$$
with the generating functions, respectively
$$
F(n)=\frac{(-1)^n}{2}, \;\;\;\;\; F(n)=n-1, \;\;\;\;\; F(n)=B_{2k}(n-1)
$$
$$
F(n)=\frac{(-1)^n}{2k+1}\sum_{u=1}^{2k+1}(2^u-1)\binom{2k+1}{u} B_{u}(n-1)^{2k+1-u}
$$
\end{enumerate}
\end{example}

\begin{theorem}  \label{theor:12}
For any polynomial $f(x)$, $x\in R$
\begin{equation} 
\lim\limits_{n\rightarrow\infty}(-1)^{n}f(n)=0   \label{eq:13}
\end{equation}
\begin{equation} 
\lim\limits_{n\rightarrow\infty}f(n)=\int\limits_{-1}^{0}f(x)dx  \label{eq:14}
\end{equation}
\end{theorem}

\section{Main Results}
\label{S4}\vspace{-4pt}

Using (\ref{eq:13}) we immediately obtain

\begin{theorem} \label{theor:13}
Let $\alpha(x)$ and $\beta(x)$ be elementary functions defined on $Z$ and satisfying the condition $\alpha(x)-\beta(x)=f(x)$, where $f(x)$ is a polynomial. Suppose that $\mu(x)$ is a function such that $\mu(x)=\alpha(x)$ if \, $2\mid x$ and $\mu(x)=\beta(x)$ if \, $2\nmid x$. Then
$$
\lim_{n\rightarrow\infty}\mu(n)=\frac{1}{2}\lim_{n\rightarrow\infty}\bigl(\alpha(n)+\beta(n)\bigr).
$$
\end{theorem}

From Theorems \ref{theor:12} and \ref{theor:13}, we get the following
\begin{theorem}  \label{theor:arithmetic}
Let $a_u=a_1+(u-1)d$, $d\geq 0$, is an arithmetic progression. Then
$$
1) \;\;\;\;\; \sum\limits_{u=1}^{\infty}a_u=\frac{5d-6a_1}{12}
$$
$$
2) \;\;\;\;\; \sum\limits_{u=1}^{\infty}(-1)^{u-1}a_u=\frac{2a_1-d}{4}
$$
\end{theorem}
\begin{example}
\begin{eqnarray}
	\sum_{u=1}^{\infty}1=1+1+1+...=-\frac{1}{2}, \;\;\;\;\;\;\; (d=0) \;\;\;\;\;\;\;\;\;\;\;\;\;\;\;\;\;\;\;\;\;\;\;\;\;\;\;\;\;\;  \textup{\cite{Ti}} \nonumber \\
	\sum_{u=1}^{\infty}u=1+2+3+...=-\frac{1}{12},\;\;\;\;\;\;\; (d=1) \;\;\;\;\;\;\;\;\;\;\;\;\;\;\;\;\;\;\;\;\;\;\;\;\;\;\;\;\;\;  \textup{\cite{Ti}} \nonumber \\
	\sum_{u=1}^{\infty}(2u-1)=1+3+5+...=\frac{1}{3},\;\;\;\;\;\;\; (d=2) \;\;\;\;\;\;\;\;\;\;\;\;\;\;\;\;\;\;\;\;\;\;\;\;\;\;\;\;\;\;  \textup{\cite{Ha}} \nonumber
\end{eqnarray}
etc.
\begin{eqnarray}
	\sum_{u=1}^{\infty}(-1)^{u-1}=1-1+1-1+...=\frac{1}{2},\;\;\;\;\;\;\; (d=0) \;\;\;\;\;\;\;\;\;\;\;\;\; \textup{\cite{Ha}} \nonumber \\
	\sum_{u=1}^{\infty}(-1)^{u-1}u=1-2+3-4+...=\frac{1}{4},\;\;\;\;\;\;\; (d=1) \;\;\;\;\;\;\;\;\;\;\;\;\; \textup{\cite{Ha}} \nonumber \\
	\sum_{u=1}^{\infty}(-1)^{u-1}(2u-1)=1-3+5-7+...=0,\;\;\;\;\;\;\; (d=2) \;\;\;\;\;\;\;\;\;\;\;\;\; \textup{\cite{Ha}} \nonumber
\end{eqnarray}
etc.
\end{example}

Let us now consider the Bernoulli polynomials. Bernoulli polynomials can be defined, in a simple way, by the symbolic equality $B_{n}(t)=(B+t)^n$, where 
the right hand side members should be expanded by the binomial theorem, and then each power $B^n$ should be replaced by $B_n$, and $B_n$ are the Bernoulli numbers
(analogously, the symbolic equality $B_n=(B+1)^n$ can be written for the Bernoulli numbers).

Bernoulli numbers play an important role in many topics of mathematics like analysis, number theory, differential topology, and in many other areas. These numbers were first introduced by Jacob Bernoulli (1654-1705) and appeared in \textit{Ars Conjectandi}, his famous treatise published posthumously in 1713,  when he studied the sums of powers of consecutive integers $1^k+2^k+3^k+...+n^k$. They often occur when expanding some simple functions in a power series. For instance, in the series
$$
\cot(x)=\frac{1}{x}-\frac{2^{2}B_2}{2!}x+\frac{2^{4}B_4}{4!}x^3-...-(-1)^{k}\frac{2^{2k}B_{2k}}{2k!}x^{2k-1}+... \;\;\; 0<|x|<\pi
$$
which appeared in astronomical works of J. Bernoulli.

Let
$$
B_{k}(n)=\frac{1}{k+1}\sum\limits_{u=0}^{k}\binom{k+1}{u}
B_{u}n^{k+1-u}=\sum\limits_{u=1}^{k}u^k
$$
be the Bernoulli polynomial.

We derive two well-known equalities,
without invoking the complex-analytical notions.

Since $B_{k}(n)-B_{k}(n-1)=n^k$, the function $f(x)=x^k$ is regular. Then, in view of (\ref{eq:11}), we have
\begin{equation}
B_{k}(-n)=(-1)^{k-1}B_{k}(n-1) \label{eq:15}
\end{equation}

On the one hand, according to (\ref{eq:15})
$$
\frac{1}{k}B_{k}(-1)=\frac{-1}{k(k+1)}\sum\limits_{u=1}^{k}(-1)^{u+k}\binom{k+1}{u}B_{u}=\frac{(-1)^{k-1}}{k}\sum\limits_{u=0}^{0}u^k=0
$$

On the other hand, in view of (\ref{eq:14})
$$
\lim\limits_{n\rightarrow\infty}B_{k-1}(n)=\int\limits_{-1}^{0}B_{k-1}(x)dx=\frac{1}{k(k+1)}\sum\limits_{u=0}^{k-1}(-1)^{u+k}\binom{k+1}{u}B_{u}=\sum\limits_{u=1}^{\infty}u^{k-1}
$$

Adding these two, we get
\begin{equation}
\sum\limits_{u=1}^{\infty}u^{k-1}=-\frac{B_k}{k} \;\;\;\;\;\;\; \forall k  \label{eq:bernoulli-1}  \;\;\;\;\;\;\;\;\;\;\;\;\;\;\;\;\;\;\;\;\;\;\;\;\;\;\;\;\;\;\;\;\;\;\;\;\;\;\;\;\;\;\;\;\;\;\;\;\;\;\;\;\;\;\;\;\; \textup{\cite{Ti}}
\end{equation}

Now taking the formula
$$
\frac{(2^{k}-1)}{k}B_{k}-\frac{(-1)^{n}}{k}\sum\limits_{u=1}^{k}(2^{u}-1)\binom{k}{u}B_{u}n^{k-u}=\sum\limits_{u=1}^{n}(-1)^{u-1}u^{k-1} \;\;\;\;\; \forall k
$$
which holds for any natural number $n$, and then passing to the limit (letting $n\rightarrow\infty$) and taking into account (\ref{eq:13}), we obtain
\begin{equation}
\sum\limits_{u=1}^{\infty}(-1)^{u-1}u^{k-1}=\frac{(2^{k}-1)}{k}B_{k} \;\;\;\;\;\;\; \forall k  \label{eq:bernoulli-2} \;\;\;\;\;\;\;\;\;\;\;\;\;\;\;\;\;\;\;\;\;\;\;\;\;\;\;\;\;\;\;\;\;\;\;\; \textup{\cite{Ha}}
\end{equation}

So, we have obtained the equalities (\ref{eq:bernoulli-1}) \cite{Ti} and (\ref{eq:bernoulli-2}) \cite{Ha}, and all the above propositions solely within
the framework of real analysis.

Combining (\ref{eq:bernoulli-1}), (\ref{eq:bernoulli-2}), Theorems \ref{theor:infsum} and \ref{theor:arithmetic}, and putting $k-1=m$, we elementarily and immediately arrive at
\begin{proposition}
For each non-negative integer $m$, we have \\
(i) for Riemann zeta function
$$
1^m+2^m+3^m+4^m+...=\zeta(-m)=-\frac{B_{m+1}}{m+1}
$$
(ii) for Dirichlet eta function
$$
1^m-2^m+3^m-4^m+...=\widetilde{\zeta}(-m)=\frac{2^{m+1}-1}{m+1}B_{m}
$$
\end{proposition}

Using (\ref{eq:15}), we can also show that all odd Bernoulli numbers $B_{k}$ (except for $B_{1}=1/2$) are equal to zero. Indeed, 
$$
B_{k}(n)-(-1)^{k-1}B_{k}(-n)=\frac{2}{k+1}\sum\limits_{u=0}^{[(k-1)/2]}\binom{k+1}{2u+1}B_{2u+1}n^{k-2u}=n^k
$$
whence we immediately obtain $B_{1}=1/2$ and $B_{2u-1}=0$, $u=1,2,...$. 

\section{Conclusion}
In conclusion, we would like to emphasize that, unlike of the many other methods that make use of either
contour integration or analytic continuation of $\zeta(s)$, our
approach, elaborated within the new theoretical direction \cite{VBT}, has allowed us to obtain the values of
the Riemann zeta function (and also the Dirichlet eta function) at non-positive integers in
a purely real way, without any notions of complex analysis and analytic continuation, and
it follows from one single concept. 
It would be interesting to apply the similar reasonings and techniques, within the setting of \cite{VBT}, to the some zeta related functions and other special functions.
We plan to study this in our subsequent works.

\vskip 1.5cm

\textbf{{Acknowledgments}}. 
I'm remembering with a deepest gratitude my teacher and supervisor Prof. Rom R. Varshamov, to whom I'm much grateful for involving me in the new research direction, 
within the framework of which this work was done. 

It is a pleasure to thank Jacques G\'elinas for bringing Sondow's paper \cite{son2} to my attention and 
informing me of \cite{landau,noah,widder} in connection with the zeroes of $\eta(s)$, as well as for useful suggestions 
on an earlier version of this article.

\end{document}